\documentclass[11pt]{article}
\usepackage{amsfonts,amscd,amssymb, amsmath}
\newtheorem{definition}{Definition}[section]
\newtheorem{lemma}[definition]{Lemma}
\newtheorem{proposition}[definition]{Proposition}
\newtheorem{corollary}[definition]{Corollary}
\newtheorem{remark}[definition]{Remark}

\newtheorem{theorem}[definition]{Theorem}
\newtheorem{example}[definition]{Example}

\def\ov{\overline}

\def\tot{\tilde{\otimes }}
\def\hot{\hat{\otimes}}
\def\oot{\ov{\otimes }}

\def\rawo\lonra{\longrightarrow}

\def\ot{\otimes}

\allowdisplaybreaks[4]

\newenvironment{proof}{{\it Proof.}}{\hfill $ \square $ \vskip 4mm}
\begin{document}
\title{Quasi-elementary $H$-Azumaya algebras arising from generalized   
(anti) Yetter-Drinfeld modules}
\author{Florin Panaite\thanks {Research
carried out while the first author was visiting the University of Antwerp,   
supported by a postdoctoral fellowship  
offered by FWO (Flemish Scientific Research Foundation). This author was 
also partially supported by the programme CEEX of the Romanian 
Ministry of Education and
Research, contract nr. 2-CEx06-11-20/2006.}\\
Institute of Mathematics of the 
Romanian Academy\\ 
PO-Box 1-764, RO-014700 Bucharest, Romania\\
e-mail: Florin.Panaite@imar.ro\\
\and 
Freddy Van Oystaeyen\\
Department of Mathematics and Computer Science\\
University of Antwerp, Middelheimlaan 1\\
B-2020 Antwerp, Belgium\\
e-mail: Francine.Schoeters@ua.ac.be}
\date{}
\maketitle
\begin{abstract}
Let $H$ be a Hopf algebra with bijective antipode, $\alpha , \beta \in 
Aut_{Hopf}(H)$ and $M$ a finite dimensional 
$(\alpha , \beta )$-Yetter-Drinfeld module. We prove that $End(M)$ 
endowed with certain structures becomes an $H$-Azumaya algebra, and the 
set of $H$-Azumaya algebras of this type is a subgroup of $BQ(k, H)$, the 
Brauer group of $H$.  
\end{abstract}
\section*{Introduction}
${\;\;\;\;}$
Let $H$ be a Hopf algebra with bijective antipode $S$ and 
$\alpha , \beta $$\; \in Aut_{Hopf}(H)$. An     
$(\alpha , \beta )$-{\it Yetter-Drinfeld  
module}, as introduced in \cite{ps}, is 
a left $H$-module right $H$-comodule $M$ with the  
following compatibility condition: 
\begin{eqnarray*}
&&(h\cdot m)_{(0)}\ot (h\cdot m)_{(1)}=h_2\cdot m_{(0)}\ot 
\beta (h_3)m_{(1)}\alpha (S^{-1}(h_1)). 
\end{eqnarray*}

This concept is a generalization of three kinds of objects appearing in the  
literature. Namely, for $\alpha =\beta =id_H$, one obtains the usual  
Yetter-Drinfeld modules; for $\alpha =S^2$, $\beta =id_H$, one obtains the  
so-called anti-Yetter-Drinfeld modules, introduced in \cite{hajac1}, 
\cite{hajac2}, \cite{jara} as coefficients for the cyclic cohomology of 
Hopf algebras defined by Connes and Moscovici in \cite{connes1},  
\cite{connes2}; finally, an $(id_H, \beta )$-Yetter-Drinfeld module is a 
generalization of a certain object $H_{\beta }$ defined in \cite{cvoz}. 
The main result in \cite{ps} is that, if we denote by ${\cal YD}(H)$ 
the disjoint union of the categories $_H{\cal YD}^H(\alpha , \beta)$ of 
$(\alpha , \beta )$-Yetter-Drinfeld modules, for all $\alpha , \beta \in 
Aut_{Hopf}(H)$, then ${\cal YD}(H)$ acquires the structure of a braided 
T-category (a concept introduced by Turaev in \cite{turaev}) over a 
certain group $G$, a semidirect product between two copies of 
$Aut_{Hopf}(H)$. Moreover, the subcategory ${\cal YD}(H)_{fd}$ 
consisting of finite dimensional objects has left and right dualities. 

The Brauer group $BQ(k, H)$ of the Hopf algebra $H$ was introduced in 
\cite{cvoz1}, by taking equivalence classes of so-called $H$-Azumaya 
algebras in the braided category $_H{\cal YD}^H$ of Yetter-Drinfeld 
modules over $H$, and using the braided product inside this 
category to define the multiplication of the group. If 
$M\in $$\;_H{\cal YD}^H$ is a finite dimensional object, 
then $End(M)$ is an  
$H$-Azumaya algebra, representing the unit element in $BQ(k, H)$. 
Also, if $H$ is finite dimensional and $\beta \in Aut_{Hopf}(H)$, the 
object $H_{\beta }$ mentioned before is not an object in $_H{\cal YD}^H$ but 
nevertheless $End(H_{\beta })$ with certain structures becomes an 
$H$-Azumaya algebra, and moreover the map 
$\beta \mapsto End(H_{\beta })$ gives a group anti-homomorphism from 
$Aut_{Hopf}(H)$ to $BQ(k, H)$, see \cite{cvoz}. 

The aim of this paper is to construct a new class of examples of 
$H$-Azumaya algebras, containing the two classes mentioned above as 
particular cases. Namely, we prove that if $\alpha , \beta\in 
Aut_{Hopf}(H)$ and $M\in $$\;_H{\cal YD}^H(\alpha , \beta )$ is finite  
dimensional, then $End(M)$ endowed with certain structures becomes an 
$H$-Azumaya algebra. The proof is rather technical and relies heavily 
on the fact that ${\cal YD}(H)_{fd}$ is a braided T-category with dualities. 
We also prove that, if we denote by $BA(k, H)$ the subset of $BQ(k, H)$ 
consisting of $H$-Azumaya algebras that can be represented as 
$End(M)$, with $M\in $$\;_H{\cal YD}^H(\alpha , \beta )$ finite 
dimensional, for some $\alpha , \beta\in Aut_{Hopf}(H)$, then $BA(k, H)$ is 
a subgroup of $BQ(k, H)$.  
\section{Preliminaries}\label{sec1}
\setcounter{equation}{0}
${\;\;\;\;}$
We work over a ground field $k$. All algebras, linear spaces,  
etc. will be over $k$; unadorned $\ot $ means $\ot_k$. Unless otherwise 
stated, $H$ will denote a Hopf algebra with bijective antipode $S$. We will 
use the version of Sweedler's sigma notation: $\Delta (h)=h_1\ot h_2$. 
For unexplained concepts and notation about Hopf algebras we refer to 
\cite{k}, \cite{m}, \cite{sw}. By $\alpha , \beta, \gamma ...$ we will 
usually denote Hopf algebra automorphisms of $H$. If $M$ is a vector space, 
a left $H$-module (respectively right $H$-comodule) structure on $M$ will 
be usually denoted by $h\otimes m\mapsto h\cdot m$ (respectively 
$m\mapsto m_{(0)}\otimes m_{(1)}$). 

We recall now some facts from \cite{ps} about 
$(\alpha , \beta )$-Yetter-Drinfeld modules. 
\begin{definition} Let $\alpha , \beta $$\;\in Aut_{Hopf}(H)$.   
An $(\alpha , \beta )$-Yetter-Drinfeld module over $H$ 
is a vector space $M$, such that $M$ is a left $H$-module (with notation 
$h\ot m\mapsto h\cdot m$) and a right $H$-comodule (with notation 
$M\rightarrow M\ot H$, $m\mapsto m_{(0)}\ot m_{(1)}$) with the following  
compatibility condition:
\begin{eqnarray}
&&(h\cdot m)_{(0)}\ot (h\cdot m)_{(1)}=h_2\cdot m_{(0)}\ot 
\beta (h_3)m_{(1)}\alpha (S^{-1}(h_1)), \label{ab}
\end{eqnarray}
for all $h\in H$ and $m\in M$. We denote by $_H{\cal YD}^H(\alpha , \beta )$  
the category of $(\alpha , \beta )$-Yetter-Drinfeld modules, morphisms 
being the $H$-linear $H$-colinear maps.
\end{definition}
\begin{remark}
As for usual Yetter-Drinfeld modules, one can see that (\ref{ab}) is 
equivalent to 
\begin{eqnarray}
&&h_1\cdot m_{(0)}\ot \beta (h_2)m_{(1)}=(h_2\cdot m)_{(0)}\ot 
(h_2\cdot m)_{(1)}\alpha (h_1). \label{abn}
\end{eqnarray}
\end{remark}
\begin{example}
For $\alpha =\beta =id_H$, we have 
$_H{\cal YD}^H(id, id)=$$\;_H{\cal YD}^H$,   
the usual category of (left-right) Yetter-Drinfeld modules. 
For $\alpha =S^2$, $\beta =id_H$, the compatibility condition (\ref{ab}) 
becomes 
\begin{eqnarray}
&&(h\cdot m)_{(0)}\ot (h\cdot m)_{(1)}=h_2\cdot m_{(0)}\ot 
h_3m_{(1)}S(h_1), \label{ayd}
\end{eqnarray}
hence $_H{\cal YD}^H(S^2, id)$ is the category of anti-Yetter-Drinfeld 
modules defined in \cite{hajac1}, \cite{hajac2}, \cite{jara}.
\end{example} 
\begin{example} \label{cvz}
For $\beta $$\;\in Aut_{Hopf}(H)$, define $H_{\beta }$ as in \cite{cvoz}, 
that is $H_{\beta }=H$, with regular right $H$-comodule structure and 
left $H$-module structure given by $h\cdot h'=\beta (h_2)h'S^{-1}(h_1)$, for 
all $h, h'\in H$. 
It was noted in \cite{cvoz} that $H_{\beta }$ satisfies a    
certain compatibility condition, which actually says that 
$H_{\beta }\in $$\;_H{\cal YD}^H(id, \beta )$. More generally, if   
$\alpha , \beta $$\;\in Aut_{Hopf}(H)$, define $H_{\alpha , \beta }$ as 
follows: $H_{\alpha , \beta }=H$, with regular right $H$-comodule structure 
and left $H$-module structure given by 
$h\cdot h'=\beta (h_2)h'\alpha (S^{-1}(h_1))$, for $h, h'\in H$. Then one 
can check that $H_{\alpha , \beta }\in $$\;_H{\cal YD}^H(\alpha , \beta )$. 
\end{example}    
\begin{example}\label{mp1}
Let $\alpha $, $\beta$$\;\in Aut_{Hopf}(H)$ and assume that there exist 
an algebra map $f:H\rightarrow k$ and a group-like element $g\in H$ such that
\begin{eqnarray} 
&&\alpha (h)=g^{-1}f(h_1)\beta(h_2)f(S(h_3))g, \;\;\;\forall \;\;h\in H. 
\label{pi2} 
\end{eqnarray}
Then one can   
check that $k\in $$\;_H{\cal YD}^H(\alpha , \beta)$, with structures 
$\;h\cdot 1=f(h)$ and $1\mapsto 1\ot g$. More generally, if $V$ is any 
vector space, then $V\in $$\;_H{\cal YD}^H(\alpha , \beta)$, with 
structures $h\cdot v=f(h)v$ and $v\mapsto v\ot g$, for all $h\in H$ and 
$v\in V$. 
\end{example}
\begin{definition}
If $\alpha $, $\beta$$\;\in Aut_{Hopf}(H)$ such that there exist $f, g$ as in 
Example \ref{mp1}, we say that $(f, g)$ is a pair in involution  
corresponding to $(\alpha,\beta)$ 
(in analogy with the concept of modular pair  
in involution due to Connes and Moscovici)  
and the $(\alpha , \beta)$-Yetter-Drinfeld modules 
$k$ and $V$ constructed in Example \ref{mp1} are denoted by $_fk^g$ and  
respectively $_fV^g$. 
\end{definition}

For instance, if $\alpha $$\;\in Aut_{Hopf}(H)$, then  
$(\varepsilon, 1)$ is a pair in involution corresponding to 
$(\alpha , \alpha )$.
\section{Tensor products and duals}\label{sec2}
\setcounter{equation}{0}
${\;\;\;\;}$    
By \cite{hajac1}, the tensor product of a   
Yetter-Drinfeld module  
with an anti-Yetter-Drinfeld module becomes an anti-Yetter-Drinfeld module.  
We generalize this result as follows:     
\begin{proposition}
Let $\alpha , \beta , \gamma $$\;\in Aut_{Hopf}(H)$  and $M, N$ two vector 
spaces which are left $H$-modules and right $H$-comodules. \\[2mm]
(i) Endow $M\ot N$ with the left $H$-module structure $h\cdot (m\ot n)=
h_1\cdot m\ot h_2\cdot n$ and the right $H$-comodule structure 
$m\ot n\mapsto (m_{(0)}\ot n_{(0)})\ot n_{(1)}m_{(1)}$ (we call these 
structures ''of type one''). If $M\in $$\;_H{\cal YD}^H(\alpha , \beta )$ and 
$N\in $$\;_H{\cal YD}^H(\beta , \gamma )$, then 
$M\ot N\in $$\;_H{\cal YD}^H(\alpha , \gamma )$; in particular, if 
$M\in $$\;_H{\cal YD}^H(S^2 , id)$ and $N\in $$\;_H{\cal YD}^H$, then 
$M\ot N\in $$\;_H{\cal YD}^H(S^2, id)$, and if 
$M\in $$\;_H{\cal YD}^H(id, \beta )$ and 
$N\in $$\;_H{\cal YD}^H(\beta , id)$, then 
$M\ot N\in $$\;_H{\cal YD}^H$.\\[2mm]
(ii) Endow $M\ot N$ with the left $H$-module structure $h\cdot (m\ot n)= 
h_2\cdot m\ot h_1\cdot n$ and the right $H$-comodule structure  
$m\ot n\mapsto (m_{(0)}\ot n_{(0)})\ot m_{(1)}n_{(1)}$ (we call these  
structures ''of type two''). 
If $M\in $$\;_H{\cal YD}^H(\alpha , \beta )$ and   
$N\in $$\;_H{\cal YD}^H(\gamma , \alpha )$, then  
$M\ot N\in $$\;_H{\cal YD}^H(\gamma , \beta )$; in particular, if  
$M\in $$\;_H{\cal YD}^H$ and $N\in $$\;_H{\cal YD}^H(S^2, id)$, then  
$M\ot N\in $$\;_H{\cal YD}^H(S^2, id)$, and if 
$M\in $$\;_H{\cal YD}^H(\alpha , id)$ and  
$N\in $$\;_H{\cal YD}^H(id, \alpha )$, then  
$M\ot N\in $$\;_H{\cal YD}^H$.
\end{proposition}  
\begin{proof}
We include here a direct proof for (i) (while (ii) is similar and 
left to the reader), an indirect proof will appear below. We compute: \\[2mm]
${\;\;\;\;}$$(h\cdot (m\ot n))_{(0)}\ot   
(h\cdot (m\ot n))_{(1)}$
\begin{eqnarray*}
&&=(h_1\cdot m\ot h_2\cdot n)_{(0)}\ot (h_1\cdot m\ot h_2\cdot n)_{(1)}\\
&&=((h_1\cdot m)_{(0)}\ot (h_2\cdot n)_{(0)})\ot (h_2\cdot n)_{(1)} 
(h_1\cdot m)_{(1)}\\
&&=(h_{(1, 2)}\cdot m_{(0)}\ot h_{(2, 2)}\cdot n_{(0)})\ot 
\gamma (h_{(2, 3)})n_{(1)}\beta (S^{-1}(h_{(2, 1)}))\beta (h_{(1, 3)})
m_{(1)}\alpha (S^{-1}(h_{(1, 1)}))\\
&&=(h_2\cdot m_{(0)}\ot h_5\cdot n_{(0)})\ot \gamma (h_6)n_{(1)}
\beta (S^{-1}(h_4)h_3)m_{(1)}\alpha (S^{-1}(h_1))\\
&&=(h_2\cdot m_{(0)}\ot h_3\cdot n_{(0)})\ot \gamma (h_4)n_{(1)}m_{(1)}
\alpha (S^{-1}(h_1))\\
&&=h_2\cdot (m\ot n)_{(0)}\ot \gamma (h_3)(m\ot n)_{(1)}\alpha (S^{-1}(h_1)),
\end{eqnarray*}
that is $M\ot N\in $$\;_H{\cal YD}^H(\alpha , \gamma )$. 
\end{proof}

In what follows, a tensor product with structures of type one will be 
denoted by $\tot $, and one with structures of type two will be denoted by 
$\oot $. 
 
By \cite{ps}, if $M\in $$\;_H{\cal YD}^H(\alpha , \beta )$ and   
$N\in $$\;_H{\cal YD}^H(\gamma , \delta )$, then  
$M\ot N$ becomes an object in $_H{\cal YD}^H(\alpha \gamma , 
\delta \gamma ^{-1}\beta \gamma )$ with the following structures:
\begin{eqnarray*}
&&h\cdot (m\ot n)=\gamma (h_1)\cdot m\ot 
\gamma ^{-1}\beta \gamma (h_2)\cdot n, \\
&&m\ot n\mapsto (m\ot n)_{(0)}\ot (m\ot n)_{(1)}:=(m_{(0)}\ot n_{(0)})\ot  
n_{(1)}m_{(1)}.
\end{eqnarray*}

This tensor product defines, on the disjoint union ${\cal YD}(H)$  
of all categories $_H{\cal YD}^H(\alpha , \beta )$, a structure of a braided  
T-category (see \cite{ps}), 
and will be denoted by $\hot $ in what follows.

We want to see what is the relation between this tensor product and 
$\tot $. We need a generalization of a result in \cite{ps},  
which states that, if $\beta $$\;\in Aut_{Hopf}(H)$, then   
$_H{\cal YD}^H(\beta , \beta )\simeq \;_H{\cal YD}^H$.  
\begin{proposition} \label{functor}
If $\alpha , \beta , \gamma $$\;\in Aut_{Hopf}(H)$, the categories      
$_H{\cal YD}^H(\alpha \beta , \gamma \beta )$ and 
$_H{\cal YD}^H(\alpha , \gamma )$ are isomorphic. 
A pair of inverse functors $(F, G)$ is given as follows. If 
$M\in $$\;_H{\cal YD}^H(\alpha \beta , \gamma \beta )$, then  
$F(M)\in $$\;_H{\cal YD}^H(\alpha , \gamma )$, 
where $F(M)=M$ as vector space, with structures 
$\;h\rightarrow m=\beta ^{-1}(h)\cdot m$ and 
$m\mapsto m_{<0>}\ot m_{<1>}:=m_{(0)}\ot m_{(1)}$, for all $h\in H$,  
$m\in M$. 
If $N\in $$\;_H{\cal YD}^H(\alpha , \gamma )$, 
then $G(N)\in $$\;_H{\cal YD}^H 
(\alpha \beta , \gamma \beta )$,   
where $G(N)=N$ as vector space, with structures 
$\;h\rightharpoonup n=\beta (h)\cdot n$ and 
$n\mapsto n^{(0)}\ot n^{(1)}:=n_{(0)}\ot n_{(1)}$, for all $h\in H$,  
$n\in N$. 
Both $F$ and $G$ act as identities on morphisms. 
\end{proposition}
\begin{proof}
Everything follows by a direct computation.
\end{proof}
\begin{corollary}
We have isomorphisms of categories: 
\begin{eqnarray*}
&&_H{\cal YD}^H(\alpha , \beta )\simeq \;
_H{\cal YD}^H(\alpha \beta ^{-1}, id),\;\;\; 
_H{\cal YD}^H(\alpha , \beta )\simeq \;
_H{\cal YD}^H(id, \beta \alpha ^{-1}),\\ 
&&_H{\cal YD}^H(\alpha , id)\simeq \;
_H{\cal YD}^H(id, \alpha ^{-1}), \;\;\; 
_H{\cal YD}^H(id, \beta )\simeq \;
_H{\cal YD}^H(\beta ^{-1}, id), 
\end{eqnarray*}
for all $\alpha , \beta \in Aut_{Hopf}(H)$.
\end{corollary}

Let now $M\in $$\;_H{\cal YD}^H(\alpha , \beta )$ and    
$N\in $$\;_H{\cal YD}^H(\beta , \gamma )$. On the one hand, we can 
consider the tensor product $M\hot N$, which is an  
object in $_H{\cal YD}^H(\alpha \beta , \gamma \beta )$. On the other 
hand, we have the tensor product $M\tot N$, which is an object in 
$_H{\cal YD}^H(\alpha , \gamma )$. Using the above formulae, one can 
then check that we have:
\begin{proposition}
$M\tot N=F(M\hot N)$.
\end{proposition}
 
Let $M$ be a finite dimensional vector space such that $M$ is a left 
$H$-module and a right $H$-comodule. Denote by $M^{\diamond}$ the 
dual vector    
space $M^*$, endowed with the following left $H$-module and right 
$H$-comodule structures:
\begin{eqnarray*} 
&&(h\cdot f) (m)=f(S(h)\cdot m), \\    
&&f_{(0)}(m)\ot f_{(1)}=f(m_{(0)})\ot S^{-1}(m_{(1)}),
\end{eqnarray*}
for all $h\in H$, $f\in M^{\diamond}$, $m\in M$, and by 
$^{\diamond}M$ the same vector space  
$M^*$ endowed with the following left $H$-module and right 
$H$-comodule structures: 
\begin{eqnarray*} 
&&(h\cdot f) (m)=f(S^{-1}(h)\cdot m), \\    
&&f_{(0)}(m)\ot f_{(1)}=f(m_{(0)})\ot S(m_{(1)}),
\end{eqnarray*}
for all $h\in H$, $f\in $$\;^{\diamond}M$, $m\in M$    
(if $M$ would be an object in $_H{\cal YD}^H$,  
then $M^{\diamond}$ and $^{\diamond}M$ would be the left and 
right duals of $M$ in $_H{\cal YD}^H$).
\begin{proposition}
If $M$ is a finite dimensional object in $_H{\cal YD}^H(\alpha , \beta )$, 
then $M^{\diamond}$ and $^{\diamond}M$ are objects in 
$_H{\cal YD}^H(\beta , \alpha )$. 
\end{proposition}
\begin{proof}
Follows by direct computation (an alternative proof will appear below).
\end{proof}
 
Recall now from \cite{ps} that, if $M$ is a finite dimensional object in 
$_H{\cal YD}^H(\alpha , \beta )$, then $M$ has left and right duals $M^*$ 
and respectively $^*M$ in the T-category ${\cal YD}(H)$; in particular,  
$M^*$ and $^*M$ are objects in 
$_H{\cal YD}^H(\alpha^{-1} , \alpha\beta^{-1}\alpha^{-1})$, 
defined as follows:  
as vector spaces they coincide both to the dual vector space of $M$, with  
structures:
\begin{eqnarray*}
&&(h\cdot f) (m)=f((\beta^{-1}\alpha^{-1}S(h))\cdot m), \\    
&&f_{(0)}(m)\ot f_{(1)}=f(m_{(0)})\ot S^{-1}(m_{(1)}), 
\end{eqnarray*}
for $M^*$, and 
\begin{eqnarray*}
&&(h\cdot f) (m)=f((\beta^{-1}\alpha^{-1}S^{-1}(h))\cdot m), \\    
&&f_{(0)}(m)\ot f_{(1)}=f(m_{(0)})\ot S(m_{(1)}), 
\end{eqnarray*} 
for $^*M$. We are interested to see how the objects 
$M^{\diamond}$, $M^*$ and    
respectively $^{\diamond}M$, $^*M$ are related. Consider the functor $F$ 
as in Proposition \ref{functor}, but this time between the categories   
$_H{\cal YD}^H(\beta (\beta ^{-1}\alpha ^{-1}), \alpha (\beta ^{-1}
\alpha ^{-1}))$ and $_H{\cal YD}^H(\beta , \alpha )$. Then, using the 
expression for $F$ and the above formulae, one can check that we have:
\begin{proposition}\label{duals}
$M^{\diamond}=F(M^*)$ and $^{\diamond}M=F(^*M)$.  
\end{proposition}
\begin{lemma} \label{ajuta}
Let $M\in $$\;_H{\cal YD}^H(\alpha , \beta )$ and    
$N\in $$\;_H{\cal YD}^H(\gamma , \delta )$ finite dimensional. Then 
the map 
\begin{eqnarray*}
&&\Psi :N^*\hot M^*\rightarrow (M\hot N)^*, \;\;\;\Psi (n^*\otimes 
m^*)(m\otimes n):=m^*(m)n^*(n), 
\end{eqnarray*} 
is an isomorphism in $_H{\cal YD}^H(\gamma ^{-1}\alpha ^{-1}, 
\alpha \beta ^{-1}\gamma \delta ^{-1}\gamma ^{-1}\alpha ^{-1})$.  
\end{lemma}
\begin{proof}
Straightforward computation.
\end{proof}
\section{Endomorphism algebras}\label{sec3}
\setcounter{equation}{0}
${\;\;\;\;}$
Let $A$ be an algebra in $_H{\cal YD}^H$. We denote by $A^{op}$ the 
(usual) opposite algebra, with multiplication $a\bullet a'=a'a$ for all 
$a, a'\in A$, and by $\overline{A}$ the $H$-opposite algebra (the 
opposite of $A$ in the category $_H{\cal YD}^H$), which equals $A$ as 
object in $_H{\cal YD}^H$ but with multiplication 
$a*a'=a'_{(0)}(a'_{(1)}\cdot a)$, for all $a, a'\in A$.

If $A$, $B$ are algebras in $_H{\cal YD}^H$, then $A\ot B$ becomes also 
an algebra in $_H{\cal YD}^H$ with the following structures:
\begin{eqnarray*}
&&h\cdot (a\ot b)=h_1\cdot a \ot h_2\cdot b,\\
&&(a\ot b)\mapsto (a_{(0)}\ot b_{(0)})\ot b_{(1)}a_{(1)}, \\
&&(a\ot b)(a'\ot b')=aa'_{(0)}\ot (a'_{(1)}\cdot b)b'.
\end{eqnarray*}

This algebra structure on $A\ot B$ (which is just the braided tensor 
product of $A$ and $B$ in the braided category $_H{\cal YD}^H$)  
is denoted by $A\# B$ and its elements  
are denoted by $a\# b$.

We introduce now endomorphism algebras associated to  
$(\alpha , \beta )$-Yetter-Drinfeld modules. 
\begin{proposition} \label{endostr}
Let $\alpha , \beta \in Aut_{Hopf}(H)$ and  
$M\in $$\;_H{\cal YD}^H(\alpha , \beta )$ finite dimensional. Then: \\
(i) $End (M)$ becomes an algebra in $_H{\cal YD}^H$, with structures:
\begin{eqnarray*}
&&(h\cdot u)(m)=\alpha ^{-1}(h_1)\cdot u(\alpha ^{-1}(S(h_2))\cdot m), \\
&&u_{(0)}(m)\ot u_{(1)}=u(m_{(0)})_{(0)}\ot S^{-1}(m_{(1)})
u(m_{(0)})_{(1)},
\end{eqnarray*}
for all $h\in H$, $u\in End (M)$, $m\in M$. \\  
(ii) $End (M)^{op}$ becomes an algebra in $_H{\cal YD}^H$, with structures:
\begin{eqnarray*}
&&(h\cdot u)(m)=\beta ^{-1}(h_2)\cdot u(\beta ^{-1}(S^{-1}(h_1))\cdot m), \\
&&u_{(0)}(m)\ot u_{(1)}=u(m_{(0)})_{(0)}\ot 
u(m_{(0)})_{(1)}S(m_{(1)}), 
\end{eqnarray*}
for all $h\in H$, $u\in End (M)^{op}$, $m\in M$. 
\end{proposition}  
\begin{proof}
Everything follows by direct computation. Note that the structures of 
$End (M)$ can be obtained in two (equivalent) ways, namely either take 
$M\hot M^*$, which is in $_H{\cal YD}^H$, and then transfer its structures  
to $End (M)$, or by taking first $M\tot M^{\diamond}$, which is in 
$_H{\cal YD}^H(\alpha , \alpha )$, transforming this into an object in 
$_H{\cal YD}^H$ via the isomorphism $_H{\cal YD}^H(\alpha , \alpha )\simeq 
\;_H{\cal YD}^H$, and finally transferring the structures to $End(M)$.  
Similarly, the structures of $End(M)^{op}$ can be obtained either by 
transferring the structures from $^*M\hot M$, which is in $_H{\cal YD}^H$,  
or by taking first $^{\diamond}M\tot M$, which is in 
$_H{\cal YD}^H(\beta , \beta )$, transforming this into an object in 
$_H{\cal YD}^H$ via the isomorphism $_H{\cal YD}^H(\beta , \beta )\simeq 
\;_H{\cal YD}^H$ and finally transferring the structures to $End(M)^{op}$.
\end{proof}
\begin{remark}
Assume that there exists a pair in involution $(f, g)$ corresponding to 
$(\alpha , \beta )$ and consider the $(\alpha , \beta )$-Yetter-Drinfeld 
module  
$_fk^g$ as in the Preliminaries. Then one can easily check that 
$End (_fk^g)$ coincides, as an algebra in $_H{\cal YD}^H$, with $k$ with  
trivial Yetter-Drinfeld structures.
\end{remark}

Let $\alpha , \beta , \gamma \in Aut_{Hopf}(H)$ and the functor $F$ as in 
Proposition \ref{functor}. Then one can easily check that we have:
\begin{corollary}\label{coco}
If $M$ is a finite dimensional object in $_H{\cal YD}^H(\alpha \beta , 
\gamma \beta )$ and consider the object $F(M)\in \;_H{\cal YD}^H(\alpha , 
\gamma )$, then $End (M)=End (F(M))$ and $End (M)^{op}=End (F(M))^{op}$ as 
algebras in $_H{\cal YD}^H$. 
\end{corollary}
\begin{corollary}
If $M\in $$\;_H{\cal YD}^H(\alpha , \beta )$ is finite dimensional,   
then 
\begin{eqnarray*}
&&End (M^*)=End (M^{\diamond}), \;\;End (^*M)=End (^{\diamond}M),\\ 
&&End (M^*)^{op}=End (M^{\diamond})^{op}, \;\;End (^*M)^{op}=
End (^{\diamond}M)^{op},
\end{eqnarray*}
as algebras in $_H{\cal YD}^H$.
\end{corollary}
\begin{proof}
Follows from Proposition \ref{duals} and Corollary \ref{coco}.
\end{proof}
\begin{corollary}\label{cocu}
Let $\alpha , \beta \in Aut_{Hopf}(H)$ and assume 
that $H$ is moreover finite   
dimensional. Then $End (H_{\alpha , \beta })=End (H_{\beta \alpha ^{-1}})$ 
as algebras in $_H{\cal YD}^H$, where $H_{\alpha , \beta }$ and 
$H_{\beta \alpha ^{-1}}$ are as in Example \ref{cvz}.
\end{corollary}
\begin{proof}
Follows from Corollary \ref{coco}, using the fact that 
$H_{\alpha , \beta }$ and $H_{\beta \alpha ^{-1}}$ correspond via the 
isomorphism of categories $_H{\cal YD}^H(\alpha , \beta )\simeq \;
_H{\cal YD}^H(id, \beta \alpha ^{-1})$.
\end{proof}

Let $M$ be a finite dimensional vector space endowed with a left 
$H$-module and a right $H$-comodule structures. Consider on $End(M)$ the 
canonical left $H$-module and right $H$-comodule structures induced by the 
structures of $M$, that is 
\begin{eqnarray*}
&&(h\cdot u)(m)=h_1\cdot u(S(h_2)\cdot m), \\
&&u_{(0)}(m)\ot u_{(1)}=u(m_{(0)})_{(0)}\ot S^{-1}(m_{(1)})
u(m_{(0)})_{(1)},
\end{eqnarray*}
for all $h\in H$, $u\in End (M)$, $m\in M$. We recall the following 
concept from \cite{cvoz}: if $A$ is an algebra in $_H{\cal YD}^H$, 
then $A$ is called {\em quasi-elementary} if there exists such an $M$ 
with the property that $End(M)$ with the above structures is an 
algebra in $_H{\cal YD}^H$ which coincides with $A$ as an algebra in  
$_H{\cal YD}^H$. 
\begin{proposition} \label{quasi}
Let $\alpha , \beta \in Aut_{Hopf}(H)$ and 
$M\in $$\;_H{\cal YD}^H(\alpha , \beta )$ finite dimensional. Then 
$End(M)$ with structures as in Proposition \ref{endostr} is a 
quasi-elementary algebra in $_H{\cal YD}^H$. 
\end{proposition}  
\begin{proof}
This is obvious if $\alpha =id_H$, because of the formulae for 
the $H$-module  
and $H$-comodule structures of $End(M)$ given in Proposition  
\ref{endostr} (i). For the general case, we consider the functor  
$F:\;_H{\cal YD}^H(\alpha , \beta )\rightarrow \; 
_H{\cal YD}^H(id_H , \beta \alpha ^{-1})$ as in Proposition \ref{functor}. 
We know from Corollary \ref{coco} that $End(M)=End(F(M))$ as algebras in 
$_H{\cal YD}^H$, and since $End(F(M))$ is quasi-elementary it follows that 
so is $End(M)$. We emphasize that the $H$-module $H$-comodule object 
making $End(M)$ quasi-elementary is {\em not} $M$ itself, but $F(M)$.   
\end{proof}
 
Recall from \cite{ps} the group $G=Aut_{Hopf}(H)\times Aut_{Hopf}(H)$ with 
multiplication $(\alpha , \beta )*(\gamma , \delta )=(\alpha \gamma , 
\delta \gamma ^{-1}\beta \gamma )$. We have the obvious result: 
\begin{lemma}\label{lulu}
The map $G\rightarrow Aut_{Hopf}(H)$, $(\alpha , \beta )\mapsto 
\beta \alpha ^{-1}$ is a group anti-homomorphism.
\end{lemma}
\begin{proposition}
If $H$ is finite dimensional, the map $(\alpha , \beta )\mapsto 
End (H_{\alpha , \beta })$ defines a group homomorphism from $G$ to the  
Brauer group $BQ(k, H)$ of $H$. 
\end{proposition}
\begin{proof}
Using Corollary \ref{cocu}, the map  
$(\alpha , \beta )\mapsto End (H_{\alpha , \beta })$ is just the composition 
between the group anti-homomorphisms $G\rightarrow Aut_{Hopf}(H)$ from 
Lemma \ref{lulu} and $Aut_{Hopf}(H)\rightarrow BQ(k, H)$, 
$\alpha \mapsto End (H_{\alpha })$ from \cite{cvoz}. 
\end{proof}

Let $\beta \in Aut_{Hopf}(H)$ and $H_{\beta }$ as in Example \ref{cvz}; 
in \cite{cvoz} was defined another object, denoted by $H'_{\beta }$, 
as follows: it has the same left $H$-module structure as $H_{\beta }$, 
and right $H$-comodule structure given by $h\mapsto h_1\ot \beta ^{-1}(h_2)$. 
It was proved then that $H'_{\beta }$ satisfies a certain compatibility 
condition, which actually says that $H'_{\beta }\in $$
\;_H{\cal YD}^H(\beta ^{-1}, id)$.

If instead of $H_{\beta }$ we take an arbitrary object $M\in $$ 
\;_H{\cal YD}^H(\alpha , \beta )$, with $\alpha , \beta \in Aut_{Hopf}(H)$,  
then the above result admits several possible generalizations; we choose 
the one that will serve our next purpose, which will be to identify the  
$H$-opposite of $End(M)$ (in case $M$ is finite dimensional), 
generalizing \cite{cvoz}, Lemma 4.5 as well as \cite{cvoz1}, 
Proposition 4.2.   
\begin{proposition}
Let $\alpha , \beta \in Aut_{Hopf}(H)$ and 
$M\in $$\;_H{\cal YD}^H(\alpha , \beta )$. Define a new object $M'$ as 
follows: $M'$ coincides with $M$ as left $H$-modules, and has a right  
$H$-comodule structure given by 
\begin{eqnarray*}
&&m\mapsto m_{<0>}\ot m_{<1>}:=m_{(0)}\ot \alpha \beta ^{-1}(m_{(1)}), 
\end{eqnarray*}
where $m\mapsto m_{(0)}\ot m_{(1)}$ is the comodule structure of $M$. 
Then $M'\in $$\;_H{\cal YD}^H(\alpha \beta ^{-1}\alpha , \alpha )$. 
\end{proposition}  
\begin{proof}
We compute:
\begin{eqnarray*}
(h\cdot m)_{<0>}\ot (h\cdot m)_{<1>}&=&
(h\cdot m)_{(0)}\ot \alpha \beta ^{-1}((h\cdot m)_{(1)})\\
&=&h_2\cdot m_{(0)}\ot \alpha \beta ^{-1}(\beta (h_3)m_{(1)}
\alpha (S^{-1}(h_1)))\\
&=&h_2\cdot m_{(0)}\ot \alpha (h_3)\alpha \beta ^{-1}(m_{(1)})
\alpha \beta ^{-1}\alpha (S^{-1}(h_1))\\
&=&h_2\cdot m_{<0>}\ot \alpha (h_3)m_{<1>}\alpha \beta ^{-1}\alpha 
(S^{-1}(h_1)),
\end{eqnarray*}
that is $M'\in $$\;_H{\cal YD}^H(\alpha \beta ^{-1}\alpha , \alpha )$.
\end{proof}
\begin{proposition}\label{p1}
Let $\alpha , \beta \in Aut_{Hopf}(H)$ and $M\in \;_H{\cal YD}^H(\alpha , 
\beta )$ finite dimensional; consider also the object 
$M'\in \;_H{\cal YD}^H(\alpha \beta ^{-1}\alpha , \alpha )$ as above. 
Define the map 
\begin{eqnarray*}
&&\tau :\overline{End (M)}\rightarrow End (M')^{op},\;\;\;\;
\tau (u)(m)=u_{(0)}(\alpha ^{-1}(u_{(1)})\cdot m), 
\end{eqnarray*}
for all $u\in End (M)$ and $m\in M'$, where $u\mapsto u_{(0)}\ot 
u_{(1)}$ is the right $H$-comodule structure of $End (M)$. Then $\tau $ is 
an isomorphism of algebras in $_H{\cal YD}^H$. 
\end{proposition}
\begin{proof}
We first prove that $\tau $ is an algebra map. We compute:
\begin{eqnarray*}
\tau (u*v)(m)&=&\tau (v_{(0)}(v_{(1)}\cdot u))(m)\\
&=&(v_{(0)}(v_{(1)}\cdot u))_{(0)}(\alpha ^{-1}
((v_{(0)}(v_{(1)}\cdot u))_{(1)})\cdot m)\\
&=&v_{(0)(0)}(v_{(1)}\cdot u)_{(0)}(\alpha ^{-1}((v_{(1)}\cdot u)_{(1)}
v_{(0)(1)})\cdot m)\\
&=&v_{(0)}(v_{(1)_2}\cdot u)_{(0)}(\alpha ^{-1}((v_{(1)_2}\cdot u)_{(1)}
v_{(1)_1})\cdot m)\\
&=&v_{(0)}(v_{(1)_3}\cdot u_{(0)})(\alpha ^{-1}(v_{(1)_4}u_{(1)}
S^{-1}(v_{(1)_2})v_{(1)_1})\cdot m)\\
&=&v_{(0)}(v_{(1)_1}\cdot u_{(0)})(\alpha ^{-1}(v_{(1)_2}u_{(1)})\cdot m)\\
&=&v_{(0)}(\alpha ^{-1}(v_{(1)_1})\cdot u_{(0)}(\alpha ^{-1}
(S(v_{(1)_2}))\alpha ^{-1}(v_{(1)_3}u_{(1)})\cdot m))\\
&=&v_{(0)}(\alpha ^{-1}(v_{(1)})\cdot u_{(0)}(\alpha ^{-1}(u_{(1)})\cdot m))\\
&=&\tau (v)(u_{(0)}(\alpha ^{-1}(u_{(1)})\cdot m))\\
&=&\tau (v)(\tau (u)(m))\\
&=&(\tau (u)\bullet \tau (v))(m), \;\;\;q.e.d.
\end{eqnarray*}
We prove now that $\tau $ is $H$-linear. We compute:
\begin{eqnarray*}
\tau (h\cdot u)(m)&=&(h\cdot u)_{(0)}
(\alpha ^{-1}((h\cdot u)_{(1)})\cdot m)\\
&=&(h_2\cdot u_{(0)})(\alpha ^{-1}(h_3u_{(1)}S^{-1}(h_1))\cdot m)\\
&=&\alpha ^{-1}(h_2)\cdot u_{(0)}(\alpha ^{-1}(S(h_3))\alpha ^{-1}(h_4
u_{(1)}S^{-1}(h_1))\cdot m)\\
&=&\alpha ^{-1}(h_2)\cdot u_{(0)}(\alpha ^{-1}
(u_{(1)}S^{-1}(h_1))\cdot m)\\
&=&\alpha ^{-1}(h_2)\cdot u_{(0)}(\alpha ^{-1}(u_{(1)})\alpha ^{-1}
(S^{-1}(h_1))\cdot m)\\
&=&\alpha ^{-1}(h_2)\cdot \tau (u)(\alpha ^{-1}
(S^{-1}(h_1))\cdot m)\\
&=&(h\cdot \tau (u))(m), \;\;\;q.e.d.
\end{eqnarray*}
We prove now that $\tau $ is $H$-colinear. 
We have to prove that $\rho (\tau (u))=\tau (u_{(0)})\ot u_{(1)}$, 
if we denote by $\rho $ the $H$-comodule structure of $End (M')^{op}$, 
that is, if we denote $\rho (v)=v^{(0)}\ot v^{(1)}$, we have to prove 
that $\tau (u)^{(0)}(m)\ot \tau (u)^{(1)}=\tau (u_{(0)})(m)\ot u_{(1)}$ 
for all $m\in M'$. Recall that the $H$-comodule structure of $M'$ is 
given by $m\mapsto m_{<0>}\ot m_{<1>}=m_{(0)}\ot \alpha \beta ^{-1}
(m_{(1)})$. First we compute:
\begin{eqnarray*}
\tau (u_{(0)})(m)\ot u_{(1)}&=&u_{(0)(0)}(\alpha ^{-1}(u_{(0)(1)})\cdot m)
\ot u_{(1)}\\
&=&u_{(0)}(\alpha ^{-1}(u_{(1)_1})\cdot m)\ot u_{(1)_2}\\
&=&u_{(0)}(\alpha ^{-1}(u_{(1)_3})\cdot m_{(0)})\\
&&\ot \alpha \beta ^{-1}(\beta (\alpha ^{-1}(u_{(1)_4}))m_{(1)_1}
\alpha (S^{-1}(\alpha ^{-1}(u_{(1)_2})))u_{(1)_1}S(m_{(1)_2}))\\
&=&u_{(0)}((\alpha ^{-1}(u_{(1)_2})\cdot m_{(0)})_{(0)})\ot 
\alpha \beta ^{-1}((\alpha ^{-1}(u_{(1)_2})\cdot m_{(0)})_{(1)}u_{(1)_1}
S(m_{(1)})).
\end{eqnarray*}
On the other hand, note that the formula for the $H$-comodule structure of 
$End (M)$ implies 
\begin{eqnarray}
&&u(m)_{(0)}\ot u(m)_{(1)}=u_{(0)}(m_{(0)})\ot m_{(1)}u_{(1)}, \label{com}
\end{eqnarray}
for all $u\in End (M)$ and $m\in M$. Now we compute:
\begin{eqnarray*}
\tau (u)^{(0)}(m)\ot \tau (u)^{(1)}&=&\tau (u)(m_{<0>})_{<0>}\ot 
\tau (u)(m_{<0>})_{<1>}S(m_{<1>})\\
&=&\tau (u)(m_{(0)})_{(0)}\ot \alpha \beta ^{-1}(\tau (u)(m_{(0)})_{(1)}
S(m_{(1)}))\\
&=&(u_{(0)}(\alpha ^{-1}(u_{(1)})\cdot m_{(0)}))_{(0)}\\
&&\ot \alpha \beta ^{-1}((u_{(0)}(\alpha ^{-1}(u_{(1)})
\cdot m_{(0)}))_{(1)}S(m_{(1)}))\\
&\overset{(\ref{com})}{=}&u_{(0)(0)}((\alpha ^{-1}(u_{(1)})
\cdot m_{(0)})_{(0)})\\ 
&&\ot \alpha \beta ^{-1}((\alpha ^{-1}(u_{(1)})
\cdot m_{(0)})_{(1)}u_{(0)(1)}S(m_{(1)})), 
\end{eqnarray*}
so the two terms are equal. The only thing left to prove is that 
$\tau $ is bijective; define the map 
\begin{eqnarray*}
&&\tau ^{-1}:End (M')^{op}\rightarrow \overline{End (M)},\;\;\;\;
\tau ^{-1}(v)(m)=v^{(0)}(\alpha ^{-1}(S(v^{(1)}))\cdot m),
\end{eqnarray*}
for all $v\in End (M')^{op}$ and $m\in M$. From the $H$-colinearity of 
$\tau $ it follows easily that $\tau ^{-1}\tau =id$. We have not been 
able to  
prove directly that $\tau \tau ^{-1}=id$; we need to prove first   
that $\tau ^{-1}$ is also $H$-colinear, that is we have to prove that 
\begin{eqnarray*}
&&\tau ^{-1}(v)(m_{(0)})_{(0)}\ot S^{-1}(m_{(1)})\tau ^{-1}(v)(m_{(0)})_{(1)}
=\tau ^{-1}(v^{(0)})(m)\ot v^{(1)},
\end{eqnarray*}
for all $v\in End (M')^{op}$ and $m\in M$. Note first that 
\begin{eqnarray*}
&&v^{(0)}(m)\ot v^{(1)}=v(m_{(0)})_{(0)}\ot \alpha \beta ^{-1}
(v(m_{(0)})_{(1)}S(m_{(1)})),
\end{eqnarray*}
which together with (\ref{com}) imply 
\begin{eqnarray}
&&v(m)_{(0)}\ot v(m)_{(1)}=v^{(0)}(m_{(0)})\ot \beta 
\alpha ^{-1}(v^{(1)})m_{(1)}. \label{lala}
\end{eqnarray}
Now we compute:\\[2mm]
${\;\;\;\;}$$\tau ^{-1}(v)(m_{(0)})_{(0)}\ot S^{-1}(m_{(1)})
\tau ^{-1}(v)(m_{(0)})_{(1)}$
\begin{eqnarray*}
&=&(v^{(0)}(\alpha ^{-1}(S(v^{(1)}))\cdot m_{(0)}))_{(0)}\ot 
S^{-1}(m_{(1)})(v^{(0)}(\alpha ^{-1}(S(v^{(1)}))\cdot m_{(0)}))_{(1)}\\
&\overset{(\ref{lala})}{=}&v^{(0)(0)}((\alpha ^{-1}(S(v^{(1)}))\cdot 
m_{(0)})_{(0)}) \\  
&&\ot S^{-1}(m_{(1)})\beta \alpha ^{-1}(v^{(0)(1)})(\alpha ^{-1}(S(v^{(1)}))
\cdot m_{(0)})_{(1)}\\
&=&v^{(0)(0)}(\alpha ^{-1}(S(v^{(1)})_2)\cdot m_{(0)(0)})\\ 
&&\ot S^{-1}(m_{(1)})\beta \alpha ^{-1}(v^{(0)(1)})\beta \alpha ^{-1}
(S(v^{(1)})_3)m_{(0)(1)}S^{-1}(S(v^{(1)})_1)\\
&=&v^{(0)}(\alpha ^{-1}(S((v^{(1)})_3))\cdot m_{(0)})\\ 
&&\ot S^{-1}(m_{(1)_2})\beta \alpha ^{-1}((v^{(1)})_1)\beta \alpha ^{-1}
(S((v^{(1)})_2))m_{(1)_1}(v^{(1)})_4\\
&=&v^{(0)}(\alpha ^{-1}(S((v^{(1)})_1))\cdot m)\ot (v^{(1)})_2\\
&=&\tau ^{-1}(v^{(0)})(m)\ot v^{(1)},\;\;q.e.d.
\end{eqnarray*}
Now from the $H$-colinearity of $\tau ^{-1}$ it follows 
easily that $\tau \tau ^{-1}=id$, hence  
$\tau $ is indeed an isomorphism with inverse $\tau ^{-1}$. 
\end{proof}

It was proved in \cite{cvoz1}, Proposition 4.7 that, if 
$M$ is a finite dimensional  
Yetter-Drinfeld module, then $End (M)^{op}$ and $End (^*M)$ are isomorphic 
as algebras in $_H{\cal YD}^H$. We generalize this result as follows: 
\begin{proposition}\label{p2}
Let $\alpha , \beta \in Aut_{Hopf}(H)$ and $M\in $$\;_H{\cal YD}^H(\alpha , 
\beta )$ finite dimensional. Then $End (M)^{op}\simeq End (^{\diamond}M) 
(=End (^*M)$) as algebras in $_H{\cal YD}^H$.
\end{proposition}
\begin{proof}
Define the map 
\begin{eqnarray*}
&&\iota :End (M)^{op}\rightarrow End (^{\diamond}M),\;\;\iota (u)=u^*, 
\end{eqnarray*}
which is obviously an algebra isomorphism. We prove now that it is 
$H$-linear. Let $u\in End (M)^{op}$, $h\in H$, $f\in $$\;^{\diamond}M$ and   
$m\in M$. Using the various formulae given before (and remembering that 
$^{\diamond}M\in $$\;_H{\cal YD}^H(\beta , \alpha )$) we compute:
\begin{eqnarray*}
\iota (h\cdot u)(f)(m)&=&(f\circ (h\cdot u))(m)\\
&=&f((h\cdot u)(m))\\
&=&f(\beta ^{-1}(h_2)\cdot u(\beta ^{-1}(S^{-1}(h_1))\cdot m)),
\end{eqnarray*}
\begin{eqnarray*}
(h\cdot \iota (u))(f)(m)&=&(\beta ^{-1}(h_1)\cdot \iota (u)
(\beta ^{-1}(S(h_2))\cdot f))(m)\\
&=&(\beta ^{-1}(h_1)\cdot ((\beta ^{-1}(S(h_2))\cdot f)\circ u))(m)\\
&=&((\beta ^{-1}(S(h_2))\cdot f)\circ u)(\beta ^{-1}(S^{-1}(h_1))\cdot m)\\
&=&(\beta ^{-1}(S(h_2))\cdot f)(u(\beta ^{-1}(S^{-1}(h_1))\cdot m))\\
&=&f(\beta ^{-1}(h_2)\cdot u(\beta ^{-1}(S^{-1}(h_1))\cdot m)),
\end{eqnarray*}
hence the two terms are equal. 
The $H$-colinearity of $\iota $ is easy to prove and  
left to the reader.
\end{proof}

We recall now some more facts from \cite{ps}. If  
$N\in $$\;_H{\cal YD}^H(\gamma , \delta )$ and $\alpha , \beta \in 
Aut_{Hopf}(H)$,  
define the object $^{(\alpha , \beta )}N=N$ as vector space, with structures 
\begin{eqnarray*}
&&h\rightharpoonup n=\gamma ^{-1}\beta \gamma \alpha ^{-1}(h)\cdot n, \\
&&n\mapsto n_{<0>}\ot n_{<1>}:=n_{(0)}\ot \alpha \beta ^{-1}(n_{(1)}). 
\end{eqnarray*}
Then $^{(\alpha , \beta )}N\in $$\;_H{\cal YD}^H
(\alpha \gamma \alpha ^{-1}, \alpha \beta ^{-1}\delta \gamma ^{-1}\beta 
\gamma \alpha ^{-1})=$$\;_H{\cal YD}^H((\alpha , \beta )*(\gamma , \delta )* 
(\alpha , \beta )^{-1})$, where $*$ is the multiplication in the group 
$G$ recalled before. Let also  
$M\in $$\;_H{\cal YD}^H(\alpha , \beta )$ and denote by $^MN=\;
^{(\alpha , \beta )}N$; then the braiding in the T-category ${\cal YD}(H)$ 
is given by the maps 
\begin{eqnarray*}
&&c_{M, N}:M\hot N\rightarrow \;^MN\hot M,\;\;\;c_{M, N}(m\ot n)=
n_{(0)}\ot \beta ^{-1}(n_{(1)})\cdot m, 
\end{eqnarray*}
which are isomorphisms in $_H{\cal YD}^H((\alpha , \beta )*
(\gamma , \delta ))$. In particular, assume that $\alpha =\beta =id_H$, 
so $M\in $$\;_H{\cal YD}^H$; then obviously $^MN=N$ as objects in 
$_H{\cal YD}^H(\gamma , \delta )$ and we have the isomorphism in 
$_H{\cal YD}^H(\gamma , \delta )$ 
\begin{eqnarray*}
&&c_{M, N}:M\hot N\rightarrow \;N\hot M,\;\;\;c_{M, N}(m\ot n)=
n_{(0)}\ot n_{(1)}\cdot m, 
\end{eqnarray*}
with inverse $c_{M, N}^{-1}(n\otimes m)=S(n_{(1)})\cdot m\otimes n_{(0)}$.

It was proved in \cite{cvoz1}, Proposition 4.3 that, 
if $M$ and $N$ are finite dimensional  
Yetter-Drinfeld modules, then $End (M)\# End (N)\simeq End (M\ot N)$ as 
algebras in $_H{\cal YD}^H$. We generalize this result as follows:
\begin{proposition}\label{p3}
If $M\in $$\;_H{\cal YD}^H(\alpha , \beta )$ and 
$N\in $$\;_H{\cal YD}^H(\gamma , \delta )$ both finite dimensional, then 
$End (M)\# End (N)\simeq End (M\hot N)$ as algebras in 
$_H{\cal YD}^H$. 
\end{proposition} 
\begin{proof}
Define the map $\phi :End (M)\# End (N)\rightarrow End (M\hot N)$ by the 
formula  
\begin{eqnarray*}
&&\phi (u\# v)(m\ot n)=u(m_{(0)})\ot (m_{(1)}\cdot v)(n),
\end{eqnarray*}
for all $u\in End (M)$, $v\in End (N)$, $m\in M$, $n\in N$, where 
$\cdot $ is the $H$-module structure of $End(N)$ as in  
Proposition \ref{endostr} (i). As in  
\cite{cvoz1} one can prove that $\phi $ is an algebra map. We prove now that 
$\phi $ is $H$-linear. We compute:
\begin{eqnarray*}
\phi (h\cdot (u\# v))(m\ot n)&=&\phi (h_1\cdot u\# h_2\cdot v)(m\ot n)\\
&=&(h_1\cdot u)(m_{(0)})\ot (m_{(1)}h_2\cdot v)(n)\\
&=&\alpha ^{-1}(h_1)\cdot u(\alpha ^{-1}(S(h_2))\cdot m_{(0)})\\ 
&& \ot \gamma ^{-1}(m_{(1)_1}h_3)\cdot v(\gamma ^{-1}(S(m_{(1)_2}h_4))
\cdot n),
\end{eqnarray*}  
${\;\;}$$(h\cdot \phi (u\# v))(m\ot n)$
\begin{eqnarray*}
&=&\gamma ^{-1}\alpha ^{-1}(h_1)\cdot (\phi (u\# v)(\gamma ^{-1}\alpha ^{-1}
(S(h_2))\cdot (m\ot n)))\\
&=&\gamma ^{-1}\alpha ^{-1}(h_1)\cdot (\phi (u\# v)(\alpha ^{-1}(S(h_3))
\cdot m\ot \gamma ^{-1}\beta \alpha ^{-1}(S(h_2))\cdot n))\\
&=&\gamma ^{-1}\alpha ^{-1}(h_1)
\cdot (u((\alpha ^{-1}(S(h_3))\cdot m)_{(0)})\\
&&\ot ((\alpha ^{-1}(S(h_3))\cdot m)_{(1)}\cdot v)(\gamma ^{-1}\beta 
\alpha ^{-1}(S(h_2))\cdot n))\\
&=&\gamma ^{-1}\alpha ^{-1}(h_1)\cdot (u(\alpha ^{-1}(S(h_4))
\cdot m_{(0)})\\
&&\ot   
(\beta \alpha ^{-1}(S(h_3))m_{(1)}h_5\cdot v)(\gamma ^{-1}\beta 
\alpha ^{-1}(S(h_2))\cdot n))\\
&=&\gamma ^{-1}\alpha ^{-1}(h_1)\cdot (u(\alpha ^{-1}(S(h_4))\cdot m_{(0)})
\ot \gamma ^{-1}\beta \alpha ^{-1}(S(h_3)_1)\gamma ^{-1}(m_{(1)_1})
\gamma ^{-1}(h_{(5, 1)})\\
&&\cdot v(\gamma ^{-1}(S(h_{(5, 2)}))
\gamma ^{-1}(S(m_{(1)_2}))\gamma ^{-1}\beta \alpha ^{-1}(S(S(h_3)_2))
\gamma ^{-1}\beta \alpha ^{-1}(S(h_2))\cdot n))\\
&=&\gamma ^{-1}\alpha ^{-1}(h_1)\cdot (u(\alpha ^{-1}(S(h_3))\cdot m_{(0)})
\ot \gamma ^{-1}\beta \alpha ^{-1}(S(h_2))\gamma ^{-1}(m_{(1)_1})
\gamma ^{-1}(h_4)\\
&&\cdot v(\gamma ^{-1}(S(h_5))\gamma ^{-1}(S(m_{(1)_2}))
\cdot n))\\
&=&\alpha ^{-1}(h_1)\cdot u(\alpha ^{-1}(S(h_4))\cdot m_{(0)})\ot 
\gamma ^{-1}\beta \alpha ^{-1}(h_2)\gamma ^{-1}\beta \alpha ^{-1}
(S(h_3))\gamma ^{-1}(m_{(1)_1}h_5)\\
&&\cdot v(\gamma ^{-1}(S(m_{(1)_2}h_6))
\cdot n)\\
&=&\alpha ^{-1}(h_1)\cdot u(\alpha ^{-1}(S(h_2))\cdot m_{(0)})\ot 
\gamma ^{-1}(m_{(1)_1}h_3)\cdot v(\gamma ^{-1}(S(m_{(1)_2}h_4))\cdot n),
\end{eqnarray*}
and we see that the two terms are equal. 
We have to prove now that $\phi $ is $H$-colinear, that is we have to 
prove that $\phi (u\# v)_{(0)}\ot \phi (u\# v)_{(1)}=\phi (u_{(0)}\# 
v_{(0)})\ot v_{(1)}u_{(1)}$. We compute: \\[2mm]
${\;\;\;\;}$$\phi (u\# v)_{(0)}(m\ot n)\ot \phi (u\# v)_{(1)}$
\begin{eqnarray*}
&=&\phi (u\# v)(m_{(0)}\ot n_{(0)})_{(0)}
\ot S^{-1}(n_{(1)}m_{(1)})\phi (u\# v)
(m_{(0)}\ot n_{(0)})_{(1)}\\
&=&u(m_{(0)})_{(0)}\ot (m_{(1)_1}\cdot v)(n_{(0)})_{(0)}\ot 
S^{-1}(n_{(1)}m_{(1)_2})(m_{(1)_1}\cdot v)(n_{(0)})_{(1)}
u(m_{(0)})_{(1)}\\
&=&u(m_{(0)})_{(0)}\ot (\gamma ^{-1}(m_{(1)_1})\cdot v(\gamma ^{-1}
(S(m_{(1)_2}))\cdot n_{(0)}))_{(0)}\\
&&\ot S^{-1}(n_{(1)}m_{(1)_3})(\gamma ^{-1}(m_{(1)_1})\cdot v(\gamma ^{-1}
(S(m_{(1)_2}))\cdot n_{(0)}))_{(1)}u(m_{(0)})_{(1)}\\
&=&u(m_{(0)})_{(0)}\ot \gamma ^{-1}(m_{(1)_2})\cdot v(\gamma ^{-1}
(S(m_{(1)_4}))\cdot n_{(0)})_{(0)}\\
&&\ot S^{-1}(n_{(1)}m_{(1)_5})\delta \gamma ^{-1}(m_{(1)_3})
v(\gamma ^{-1}(S(m_{(1)_4}))\cdot n_{(0)})_{(1)}S^{-1}(m_{(1)_1})
u(m_{(0)})_{(1)},
\end{eqnarray*}
and on the other hand:\\[2mm]
${\;\;\;\;}$$\phi (u_{(0)}\# v_{(0)})(m\ot n)\ot v_{(1)}u_{(1)}$
\begin{eqnarray*}
&=&u_{(0)}(m_{(0)})\ot (m_{(1)}\cdot v_{(0)})(n)\ot v_{(1)}u_{(1)}\\
&=&u_{(0)}(m_{(0)})\ot \gamma ^{-1}(m_{(1)_1})\cdot v_{(0)}(\gamma ^{-1}
(S(m_{(1)_2}))\cdot n)\ot v_{(1)}u_{(1)}\\
&=&u(m_{(0)(0)})_{(0)}\ot \gamma ^{-1}(m_{(1)_1})\cdot v((\gamma ^{-1}
(S(m_{(1)_2}))\cdot n)_{(0)})_{(0)}\\
&&\ot S^{-1}((\gamma ^{-1}(S(m_{(1)_2}))\cdot n)_{(1)})v((\gamma ^{-1}
(S(m_{(1)_2}))\cdot n)_{(0)})_{(1)}S^{-1}(m_{(0)(1)})u(m_{(0)(0)})_{(1)}\\
&=&u(m_{(0)})_{(0)}\ot \gamma ^{-1}(m_{(1)_2})\cdot v(\gamma ^{-1}
(S(m_{(1)_3})_2)\cdot n_{(0)})_{(0)}\\
&&\ot S^{-1}(\delta \gamma ^{-1}(S(m_{(1)_3})_3)n_{(1)}S^{-1}
(S(m_{(1)_3})_1))v(\gamma ^{-1}(S(m_{(1)_3})_2)\cdot n_{(0)})_{(1)}\\
&&S^{-1}(m_{(1)_1})u(m_{(0)})_{(1)}\\
&=&u(m_{(0)})_{(0)}\ot \gamma ^{-1}(m_{(1)_2})\cdot v(\gamma ^{-1}
(S(m_{(1)_4}))\cdot n_{(0)})_{(0)}\\
&&\ot S^{-1}(n_{(1)}m_{(1)_5})\delta \gamma ^{-1}(m_{(1)_3})
v(\gamma ^{-1}(S(m_{(1)_4}))\cdot n_{(0)})_{(1)}S^{-1}(m_{(1)_1})
u(m_{(0)})_{(1)},
\end{eqnarray*}
hence the two terms are equal. The only thing left to prove is that 
$\phi $ is bijective; we give a proof similar to the one in 
\cite{cvoz1}. Namely, one can check that $\phi $ coincides with the 
composition of the following isomorphisms: 
\begin{eqnarray*}
End(M)\hot End(N)&\simeq &M\hot M^*\hot N\hot N^*\\
&\simeq &M\hot N\hot N^*\hot M^*\\
&\simeq &M\hot N\hot (M\hot N)^*\\
&\simeq &End(M\hot N), 
\end{eqnarray*}
where the first and the last are the canonical linear isomorphisms, the 
second is $id_M\otimes c_{N\hot N^*, M^*}^{-1}$ and the third is 
$id_{M\hot N}\otimes \Psi $, where $\Psi $ is the isomorphism defined 
in Lemma \ref{ajuta}.
\end{proof}

Let us recall from \cite{cvoz} that, if $A$ is an algebra in $_H{\cal YD}^H$ 
and $\mu \in Aut_{Hopf}(H)$, we can define a new algebra in $_H{\cal YD}^H$, 
denoted by $A(\mu )$, which equals $A$ as an algebra, but with $H$-structures 
$(A(\mu ), \rightharpoondown , \rho ')$ given by 
$h\rightharpoondown a=\mu (h) \cdot a$ and $\rho '(a)=a_{<0>}\ot  
a_{<1>}:=a_{(0)}\ot \mu ^{-1}(a_{(1)})$, 
for all $a\in A(\mu )$, $h\in H$. 
\begin{proposition}\label{p4}
Let $N\in $$\;_H{\cal YD}^H(\gamma , \delta )$ finite dimensional 
and $\alpha , \beta \in Aut_{Hopf}(H)$. Then $End (^{(\alpha , \beta )}N)=
End (N)(\beta \alpha ^{-1})$ as algebras in $_H{\cal YD}^H$. 
\end{proposition}  
\begin{proof}
We compute the structures of $End (^{(\alpha , \beta )}N)$. Let $h\in H$, 
$u\in End (N)$, $n\in N$; we have: 
\begin{eqnarray*}
(h\cdot u)(n)&=&\alpha \gamma ^{-1}\alpha ^{-1}(h_1)\rightharpoonup 
u(\alpha \gamma ^{-1}\alpha ^{-1}(S(h_2))\rightharpoonup n)\\
&=&\gamma ^{-1}\beta \gamma \alpha ^{-1}\alpha \gamma ^{-1}\alpha ^{-1}(h_1)
\cdot u(\gamma ^{-1}\beta \gamma \alpha ^{-1}\alpha \gamma ^{-1}\alpha ^{-1}
(S(h_2))\cdot n)\\
&=&\gamma ^{-1}\beta \alpha ^{-1}(h_1)\cdot u(\gamma ^{-1}\beta \alpha ^{-1}
(S(h_2))\cdot n), 
\end{eqnarray*}
\begin{eqnarray*}
u_{[0]}(n)\ot u_{[1]}&=&u(n_{<0>})_{<0>}\ot S^{-1}(n_{<1>})u(n_{<0>})_{<1>}\\
&=&u(n_{(0)})_{(0)}\ot \alpha \beta ^{-1}(S^{-1}(n_{(1)})u(n_{(0)})_{(1)}),
\end{eqnarray*}
while the structures of $End (N)(\beta \alpha ^{-1})$ are: 
\begin{eqnarray*}
(h\rightharpoondown u)(n)
&=&(\beta \alpha ^{-1}(h)\cdot u)(n)\\
&=&\gamma ^{-1}\beta \alpha ^{-1}(h_1)\cdot u(\gamma ^{-1}\beta \alpha ^{-1}
(S(h_2))\cdot n), 
\end{eqnarray*}
\begin{eqnarray*}
u_{<0>}(n)\ot u_{<1>}&=&u_{(0)}(n)\ot \alpha \beta ^{-1}(u_{(1)})\\
&=&u(n_{(0)})_{(0)}\ot \alpha \beta ^{-1}(S^{-1}(n_{(1)})u(n_{(0)})_{(1)}),
\end{eqnarray*}
and we are done.
\end{proof}

It was proved in \cite{cvoz1} that, if $M$ and $N$ are finite dimensional 
Yetter-Drinfeld modules, then $End (M)\# End (N)\simeq End (N)\# End (M)$  
as algebras in $_H{\cal YD}^H$. By using Proposition \ref{p3} and the 
isomorphisms $c_{M, N}$ recalled above, we obtain the following 
generalization of this fact:
\begin{proposition}\label{p5}
Let $M\in $$\;_H{\cal YD}^H(\alpha , \beta )$ and $N\in $$\;
_H{\cal YD}^H(\gamma , \delta )$, both finite dimensional. Then  
$End (M)\# End (N)\simeq End (^M N)\# End (M)$ as algebras in 
$_H{\cal YD}^H$.
\end{proposition}
\section{$H$-Azumaya algebras and a subgroup of the Brauer group}\label{sec4}
\setcounter{equation}{0}
${\;\;\;\;}$
We begin by recalling several facts from \cite{cvoz1} about $H$-Azumaya 
algebras and the Brauer group of a Hopf algebra $H$. 

Let $A$ be a finite dimensional algebra in $_H{\cal YD}^H$ and consider 
the maps
\begin{eqnarray*}
&&F:A\# \overline{A}\rightarrow End(A), \;\;\;F(a\# b)(c)=ac_{(0)}
(c_{(1)}\cdot b), \\
&&G:\overline{A}\# A\rightarrow End(A)^{op}, \;\;\;
G(a\# b)(c)=a_{(0)}(a_{(1)}\cdot c)b,
\end{eqnarray*} 
for all $a, b, c\in A$, which are algebra maps in $_H{\cal YD}^H$.  
In case $F$ and $G$ are bijective, $A$ is called $H$-Azumaya. 
If $M$ is a finite dimensional object in $_H{\cal YD}^H$, then 
$End(M)$ is an $H$-Azumaya algebra. If $A$ and $B$ are $H$-Azumaya, 
then so are $A\# B$ and $\overline{A}$. Two $H$-Azumaya algebras 
$A$ and $B$ are called Brauer equivalent (and denoted $A\sim B$) if 
there exist $M, N\in $$\;_H{\cal YD}^H$ finite dimensional such that 
$A\# End(M)\simeq B\# End(N)$ as algebras in $_H{\cal YD}^H$. The 
relation $\sim $ is an equivalence relation which respects the operation 
$\# $. The quotient set is a group with multiplication induced by $\# $ and 
inverse induced by $A\mapsto \overline{A}$. This group is denoted by 
$BQ(k, H)$ and called the Brauer group of $H$. The class of an 
$H$-Azumaya algebra $A$ in $BQ(k, H)$ is denoted by $[A]$. 

We have now all the necessary ingredients to prove the main result 
of this paper:   
\begin{theorem}
Let $\alpha , \beta \in Aut_{Hopf}(H)$ and $M\in $$\;_H{\cal YD}^H(\alpha , 
\beta )$ a finite dimensional object. Then $End (M)$, with 
structures as in Proposition \ref{endostr}, is an $H$-Azumaya algebra. 
\end{theorem}
\begin{proof}
We prove that the map 
\begin{eqnarray*}
&&F:End(M)\# \overline{End(M)}\rightarrow End(End(M)),\\ 
&&F(a\# b)(c)=ac_{(0)}(c_{(1)}\cdot b),\;\;\;\forall \;\;a, b, c\in End(M),  
\end{eqnarray*} 
is bijective (the proof that the other map $G$ is bijective is similar 
and left to the reader). By Propositions \ref{p1}, \ref{p2} and  
\ref{p3}, we obtain that 
$End(M)\# \overline{End(M)}\simeq End(M\hot \;^{\diamond}(M'))$ as   
algebras in $_H{\cal YD}^H$; in particular, it follows that 
$End(M)\# \overline{End(M)}$ is a simple ring, and since $F$ is an 
algebra map it follows that $F$ is injective, and hence bijective, as 
$dim_k(End(M)\# \overline{End(M)})=dim_k(End(End(M)))$, finishing  
the proof.  
\end{proof}
\begin{corollary} We denote by $BA(k, H)$ the subset of $BQ(k, H)$ 
consisting of $H$-Azumaya algebras that can be represented as $End (M)$, 
with $M\in $$\;_H{\cal YD}^H(\alpha , \beta )$ finite dimensional, for some  
$\alpha , \beta \in Aut_{Hopf}(H)$. Then $BA(k, H)$ is a subgroup of 
$BQ(k, H)$. Moreover, if $H$ is finite dimensional, the image of the 
group anti-homomorphism from \cite{cvoz}, 
$Aut_{Hopf}(H)\rightarrow BQ(k, H)$, $\alpha \mapsto [End(H_{\alpha })]$, 
is contained in $BA(k, H)$. 
\end{corollary}
\begin{proof}
Follows immediately by using Propositions \ref{p1}, \ref{p2} and  
\ref{p3}. 
\end{proof}
\begin{remark}
Following \cite{cvoz}, we denote $BT(k, H)$ the subset of $BQ(k, H)$ 
consisting of classes that are represented by quasi-elementary 
$H$-Azumaya algebras. It was noted in \cite{cvoz} that $BT(k, H)$ is 
closed under multiplication, but it is not known whether it is a  
subgroup of $BQ(k, H)$. By Proposition \ref{quasi},  
$BA(k, H)\subseteq BT(k, H)$. Thus, by considering only those 
quasi-elementary $H$-Azumaya algebras that are represented as 
$End(M)$, with $M$ a finite dimensional object in some 
$_H{\cal YD}^H(\alpha , \beta )$, we do obtain a subgroup of $BQ(k, H)$.    
\end{remark}

We recall from \cite{cvoz} that the construction $A\mapsto A(\mu )$ 
recalled before defines a group action of $Aut_{Hopf}(H)$ on 
$BQ(k, H)$, by $\mu ([A])=[A(\mu )]$ for $\mu \in Aut_{Hopf}(H)$ and 
$[A]\in BQ(k, H)$. As a consequence of Proposition \ref{p4}, we obtain: 
\begin{corollary}
The above action induces a group action of $Aut_{Hopf}(H)$ on $BA(k, H)$. 
\end{corollary}
\begin{proposition}
Assume that there exists $(f, g)$ a pair in involution corresponding to 
$(\alpha , \beta )$ and let $M\in $$\;_H{\cal YD}^H(\alpha , \beta )$  
finite dimensional. Then $[End (M)]=1$ in the Brauer group. 
\end{proposition}
\begin{proof}
By \cite{ps}, Theorem 5.1, $M$ is isomorphic, as  
$(\alpha , \beta )$-Yetter-Drinfeld modules, with $_fk^g\hot N$, where 
$N\in $$\;_H{\cal YD}^H$. Thus $End (M)\simeq End (_fk^g\hot N)
\simeq End (_fk^g)\# End (N)=k\# End (N)$, hence in the Brauer group 
we get $[End (M)]=[k][End (N)]=1\cdot 1=1$.
\end{proof}

\end{document}